\documentclass[10pt]{amsart}
\usepackage{amssymb}
\usepackage{amsmath,amssymb,amsfonts,amsthm,
latexsym, amscd, amsfonts, epsfig}
\usepackage{mathrsfs}
\usepackage{color}
\usepackage{eucal}
\usepackage{eufrak}
\usepackage[all]{xypic}
\usepackage{xspace}

\theoremstyle{plain}
\newtheorem{theorem}{Theorem}

\newtheorem{lemma}{Lemma}

\theoremstyle{definition}
\newtheorem{definition}{Definition}

\theoremstyle{example}

\theoremstyle{remark}

\numberwithin{equation}{section}

\begin{document}

\title[Central and Local Limit Theorems for RNA Structures]
      {Central and Local Limit Theorems for RNA Structures}
\author{Emma Y. Jin and Christian M. Reidys$^{\,\star}$}
\address{Center for Combinatorics, LPMC-TJKLC \\
         Nankai University  \\
         Tianjin 300071\\
         P.R.~China\\
         Phone: *86-22-2350-6800\\
         Fax:   *86-22-2350-9272}
\email{reidys@nankai.edu.cn}
\thanks{}
\date{July 2007}
\begin{abstract}
A $k$-noncrossing RNA pseudoknot structure is a graph over
$\{1,\dots,n\}$ without $1$-arcs, i.e.~arcs of the form $(i,i+1)$
and in which there exists no $k$-set of mutually intersecting arcs.
In particular, RNA secondary structures are $2$-noncrossing RNA
structures. In this paper we prove a central and a local limit
theorem for the distribution of the numbers of $3$-noncrossing RNA
structures over $n$ nucleotides with exactly $h$ bonds. We will
build on the results of \cite{Reidys:07rna1} and
\cite{Reidys:07rna2}, where the generating function of
$k$-noncrossing RNA pseudoknot structures and the asymptotics for
its coefficients have been derived. The results of this paper explain
the findings on the numbers of arcs of RNA secondary structures
obtained by molecular folding algorithms and predict the distributions
for $k$-noncrossing RNA folding algorithms which are currently being
developed.
\end{abstract}
\maketitle
{{\small
}}


\section{Introduction}

\begin{figure}[ht]
\centerline{ \epsfig{file=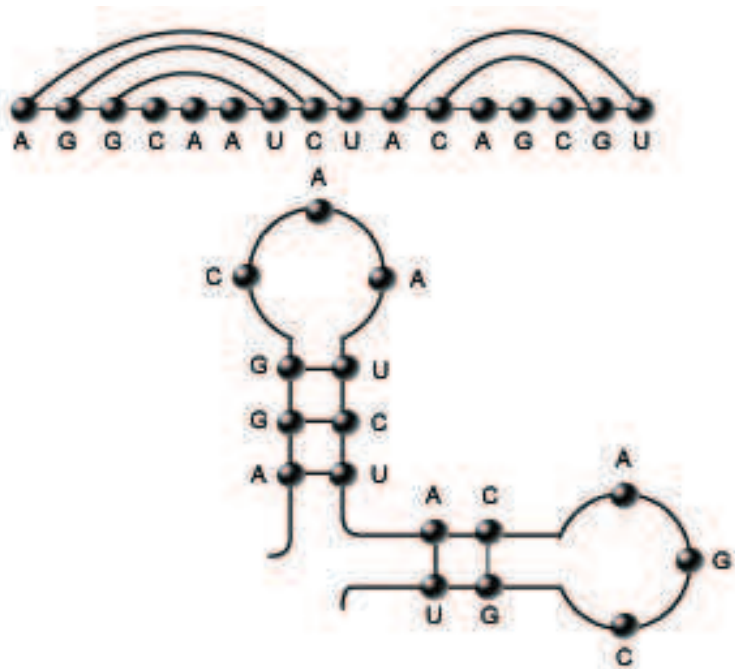,width=0.7\textwidth}\hskip15pt }
\caption{\small RNA secondary structures. Diagram representation
(top): the primary sequence, {\bf AGGCAAUCUACAGCGU}, is drawn
horizontally and its backbone bonds are ignored. All bonds are drawn
in the upper half-plane. Secondary structures have the property that
no two arcs intersect and all arcs have minimum length $2$. Outer
planar graph representation (bottom).}\label{F:1}
\end{figure}

An RNA molecule consists of the primary sequence of the four
nucleotides {\bf A}, {\bf G}, {\bf U} and {\bf C} together with the
Watson-Crick ({\bf A-U}, {\bf G-C}) and ({\bf U-G}) base pairing
rules. The latter specify the pairs of nucleotides that can
potentially form bonds. Single stranded RNA molecules form helical
structures whose bonds satisfy the above base pairing rules and
which, in many cases, determine their function. For instance RNA
ribosomes are capable of catalytic activity, cleaving other RNA
molecules. Not all possible bonds are realized, though. Due to
bio-physical constraints and the chemistry of Watson-Crick base
pairs there exist rather severe constraints on the bonds of an RNA
molecule. In light of this three decades ago Waterman {\it et.al.}
pioneered the concept of RNA secondary structures
\cite{Waterman:78a,Waterman:94a}, being
subject to the most strict combinatorial constraints. Any structure
can be represented by drawing the primary sequence horizontally,
ignoring all chemical bonds of its
backbone, see Fig.~\ref{F:1}. Then one draws all bonds, satisfying
the Watson-Crick base pairing rules as arcs in the upper half-plane,
effectively identifying structure with the set of all arcs. In this
representation, RNA secondary structures have no $1$-arcs, i.e.~arcs
of the form $(i,i+1)$ and no two arcs $(i_1,j_1)$, $(i_2,j_2)$,
where $i_1<j_1$ and $i_2<j_2$ with the property $i_1<i_2<j_1<j_2$.
In other words there exist no two arcs that cross in the diagram
representation of the structure.
It is well-known that there exist additional types of nucleotide
interactions \cite{Science:05a}. These bonds are called pseudoknots
\cite{Westhof:92a} and occur in functional RNA (RNAseP
\cite{Loria:96a}), ribosomal RNA \cite{Konings:95a} and are
conserved in the catalytic core of group I introns. Pseudoknots
appear in plant viral RNAs pseudo-knots and in {\it in vitro} RNA
evolution \cite{Tuerk:92} experiments have produced families of RNA
structures with pseudoknot motifs, when binding HIV-1 reverse
transcriptase. Important mechanisms like ribosomal frame shifting
\cite{Chamorro:91a} also involve pseudoknot interactions.
$k$-noncrossing RNA structures introduced in \cite{Reidys:07rna1}
capture these pseudoknot bonds and generalize the concept of the RNA
secondary structures in a natural way. In the diagram representation
$k$-noncrossing RNA structure has no $1$-arcs and contains at most
$k-1$ mutually crossing arcs.
\begin{figure}[ht]
\centerline{ \epsfig{file=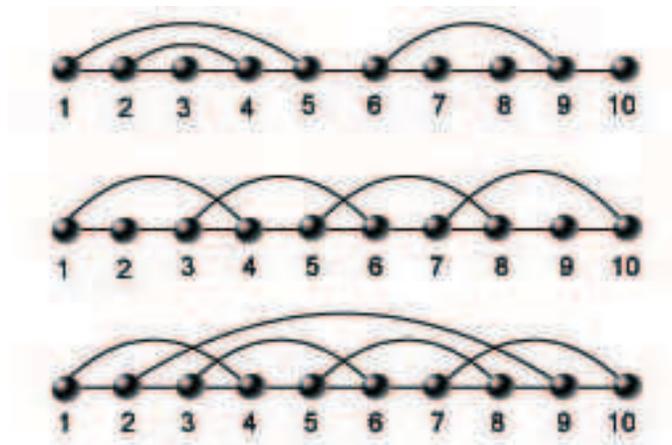,width=0.6\textwidth}\hskip15pt }
\caption{\small $k$-noncrossing RNA structures. (a) secondary
structure (with isolated labels $3,7,8,10$), (b) planar
$3$-noncrossing RNA structure, $2,9$ being isolated (c) the smallest
non-planar $3$-noncrossing structure } \label{F:2}
\end{figure}

\begin{figure}[ht]
\centerline{ \epsfig{file=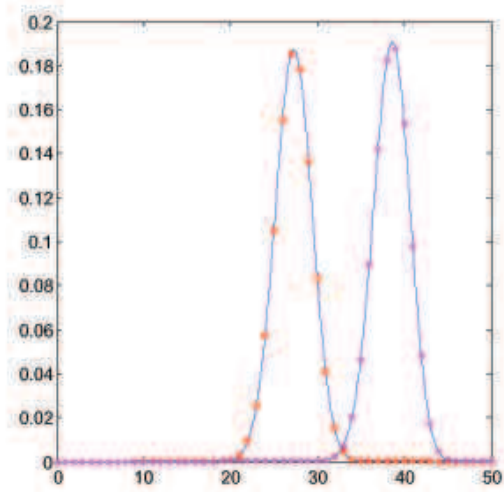,width=0.5\textwidth}\hskip15pt }
\caption{\small Central limit theorem of $2$-noncrossing and
$3$-noncrossing RNA structures: both random variables are normalized
to ${\sf{S}}'_{2}(n,h)/{\sf{S}}_{2}(n)$and
${\sf{S}}'_{3}(n,h)/{\sf{S}}_{3}(n)$, respectively. In case of
$n=100$, for $2$-noncrossing RNA structures we have a mean of
$0.276393\,n=27.6393$ and variance $0.044721\, n=4.4721$ (left
curve), while for $3$-noncrossing RNA structures mean $0.39089\,
n=39.089$ and variance $0.041565\, n=4.1565$ (right curve). The red
dots and magenta dots represent the values
${\sf{S}}'_{2}(n,h)/{\sf{S}}_{2}(n)$ and
${\sf{S}}'_{3}(n,h)/{\sf{S}}_{3}(n)$, respectively. } \label{F:4}
\end{figure}

The starting point of this paper was the experimental finding that
$3$-noncrossing RNA structures for random sequences of length $100$
over the nucleotides {\bf A}, {\bf G}, {\bf U} and {\bf C} exhibited
sharply concentrated numbers of arcs (centered at $39$). It was furthermore
intriguing that the numbers of arcs were significantly higher than
those in RNA secondary structures. While it is evident that
$3$-noncrossing RNA structures have more arcs than secondary
structures, the jump from $27$ to $39$ (for $n=100$ ) with a maximum
number of $50$ arcs was not anticipated. Since all these quantities
were via the generating functions for $k$-noncrossing RNA structures
in \cite{Reidys:07rna1} explicitly known we could easily confirm that
the numbers of $3$-noncrossing RNA structures with exactly $h$ arcs,
${\sf S}'_3 (n,h)$ satisfy indeed almost ``perfectly'' a Gaussian
distribution with a mean of $39$, see Fig.~\ref{F:3}. We also found
that a central limit theorem holds for RNA secondary structures with
$h$ arcs, see Figure~\ref{F:4}. These observation motivated us to
understand how and why these limit distributions arise, which is
what the present paper is about. Our main results can be summarized
as follows:
\begin{figure}[ht]
\centerline{ \epsfig{file=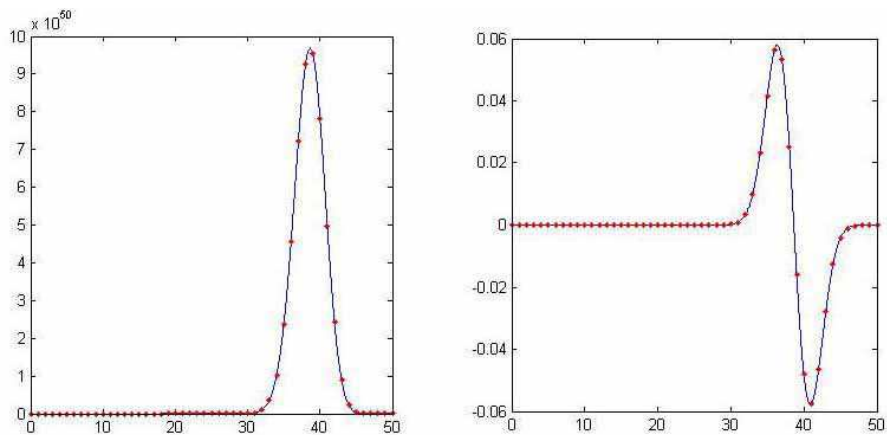,width=0.8\textwidth}}
\caption{\small Central limit theorem and local limit theorem for
3-noncrossing RNA structures of length $n=100$ with exactly $h$
arcs: we display the central limit theorem (left) for ${\sf S}'_{3}(100,h),
h=1,2,\cdots 50$ (labeled by red dots) with mean
$0.39089\cdot 100=39.089$ and variance $0.041565\cdot 100=4.1565$,
and for the local limit theorem (right), we display the difference
$\sqrt{4.1565}\ \mathbb{P}\left(\frac{X_{n}-39.089}{\sqrt{4.1565}}
=x\right)-\frac{1}{\sqrt{2 \pi}}e^{-\frac{x^2}{2}}$ which is maximal
close to the peak of the distribution.} \label{F:3}
\end{figure}

{\bf Theorem.}
{\it Let ${\sf S}_{3}'(n,h)$ denote the number of $3$-noncrossing RNA
structures with exactly $h$ arcs. Then the random variable $X_{n}$
having distribution $\mathbb{P}(X_{n}=h)={\sf S}_{3}'(n,h)/{\sf S}_{3}(n)$
satisfies a central and local limit theorem with mean $0.39089\,n$
and variance $0.041565\,n$.
}

Our particular strategy is rooted in our recent work on asymptotic
enumeration of $k$-noncrossing RNA structures \cite{Reidys:07rna2}
and a paper of Bender \cite{Bender:73} who showed how such central
limit theorems arise in case of singularities that are poles.
In order to put our results into context let us provide some
background on central and local limit theorems.
Suppose we are given a set $A_{n}$ (of size $a_n$).
For instance let $A_{n}$ be the set of subsets of $\{1,\dots,n\}$.
Suppose further we are given $A_{n,k}$ (of size $a_{n,k}$),
$k \in \mathbb{N}$ representing
a disjoint set partition of $A_{n}$. For instance let $A_{n,k}$ be
the number of subsets with exactly $k$ elements. Consider the random
variable $\xi_{n}$ having the probability distribution
$\mathbb{P}(\xi_{n}=k)=a_{n,k}/a_{n}$, then the corresponding
probability generating function is given by
$$
\sum_{k \ge 0}\mathbb{P}(\xi_{n}=k)w^k=\sum_{k \ge
0}\frac{a_{n,k}}{a_{n}}w^k=\frac{\sum_{k \ge 0}a_{n,k}w^k}{\sum_{k
\ge 0}a_{n,k}1^k} \ .
$$
Let $\varphi_{n}(w)=\sum_{k \ge 0}a_{n,k}w^k$, then
$\frac{\varphi_{n}(w)} {\varphi_{n}(1)}$ is the probability
generating function of $\xi_{n}$ and
$$
f(z,w)=\sum_{n \ge 0}\varphi_{n}(w)z^n
=\sum_{n \ge 0}\sum_{k\ge 0}a_{n,k}w^kz^n
$$
is called the bivariate generating function.
For instance, in our example we have $\mathbb{P}(\xi_{n})=\frac{{n \choose
k}}{2^n}$ and the resulting bivariate generating function is
\begin{equation}\label{E:bi}
\sum_{n \ge 0}\sum_{k \le n}{n \choose k}w^kz^n=\frac{1}{1-z(1+w)} \ .
\end{equation}
The key idea consists in considering $f(z,w)$ as being parameterized by
$w$ and to study the change of its singularity in an $\epsilon$-disc
centered at $w=1$. Indeed the moment generating function is given by
$$
E(e^{s\xi_{n}})=\sum_{k
\ge 0}\frac{a_{n,k}}{a_{n}}e^{sk}=\frac{\varphi_{n}(e^s)}{\varphi_{n}(1)}
=\frac{[z^n]f(z,e^s)}{[z^n]f(z,1)}
$$
and $\frac{[z^n]f(z,e^{it})}{[z^n]f(z,1)}=E(e^{it\xi_{n}})$ is the
characteristic function of $\xi_{n}$. This shows that the
coefficients of $f(z,w)$ control the distribution, which can, for
large $n$, be obtained via singularity analysis. The resulting analysis
can be amazingly simple. Let us showcase this in the case of the binomial
distribution. Here we have
the bivariate generating function $\sum_{n \ge 0}\sum_{k \le n}{n
\choose k}w^kz^n=\frac{1}{1-z(1+w)}$, eq.~(\ref{E:bi}). The simple
pole $r(s)$ of $f(z,e^s)$ is $\frac{1}{1+e^s}$. Observe that
$\frac{\varphi_{n}(e^s)}{\varphi_{n}(1)}\sim (\frac{r(0)}{r(s)})^n$
holds for $s$ uniformly in a neighborhood of 0, and Taylor expansion
shows
$$
\frac{\varphi_{n}(e^{it})}{\varphi_{n}(1)}\sim \exp(i\cdot
\frac{n}{2}\cdot t-\frac{1}{2}\cdot \frac{n}{4}\cdot t^2+O(t^3))
$$
uniformly for $t$ for any arbitrary finite interval.
It remains to apply the L\'{e}vy-Cram\'{e}r theorem
(Theorem~\ref{T:K}) to the normalized characteristic function
of the random variable $\frac{\eta_{n}-\frac{n}{2}}{\sqrt{\frac{n}{4}}}$,
which yields the asymptotic normality of $\eta_{n}$. Thus ${n \choose k}$
is asymptotically normal distributed with mean $\frac{n}{2}$ and
variance $\frac{n}{4}$. As it turns out we will have to work a bit harder
to prove our main result. The complication is due to the fact that the
generating function for $3$-noncrossing RNA structures is much more
complex (and fascinating) than the bivariate function of eq.~(\ref{E:bi})
which has a simple pole as dominant singularity. For instance, the
singularity of the generating function for $3$-noncrossing RNA structures
is {\it not} a pole but of algebraic-logarithmic type
\cite{Flajolet:05,Gao:92,Handbook:95}.

Our two main results, Theorem~\ref{T:gauss} in Section~\ref{S:CL}
and  Theorem~\ref{T:local} in Section~\ref{S:LL} shed light on the
distribution of $3$-noncrossing RNA structures from a global and local
perspective. A central limit theorem represents the global
perspective on the limiting distribution of some random variable $X_{n}$:
$$
\lim_{n\to\infty}\mathbb{P}\left(\frac{X_{n}-\mu_{n}}{\sigma_{n}}<x\right)
=
\frac{1}{\sqrt{2\pi}}\int_{-\infty}^{x}e^{-\frac{t^2}{2}}dt \ .
$$
Bender observed in \cite{Bender:73} that a central limit theorem
combined with certain smoothness conditions on the coefficients
$a_{n,k}$ implies a local limit theorem which considers the difference
between $\mathbb{P}(x \le \frac{X_{n}-\mu_{n}} {\sigma_{n}}<x+1)$ and
$\frac{1}{\sqrt{2\pi}}\int_{x}^{x+1}e^{-\frac{t^2}{2}}dt$ as $n$ tends
to infinity.
To be precise, $X_{n}$ satisfies a local limit theorem on some set $S\subset
\mathbb{R}$ if and only if
$$
\lim_{n \to \infty}\sup_{x\in S}\left|\sigma_{n}
\mathbb{P}\left(\frac{X_{n}-\mu_{n}}{\sigma_{n}}=x\right)-\frac{1}{\sqrt{2
\pi}}e^{-\frac{x^2}{2}}\right|=0
$$
holds and we say $X_n$ satisfies a local limit theorem for some
$S=\{x\in\mathbb{R}\mid x=o(\sqrt{n})\}$.
Why is the smoothness of the $a_{n,k}$ so important?
Suppose $a_{n,k}={n \choose k}+(-1)^k2{n \choose k}$, then
it follows in analogy to our above argument that a central limit
theorem with mean $\frac{1}{4}n$ and variance $\frac{1}{4}n$ holds.
However, $\eta_{n}$ does not satisfy a local limit theorem, since
$$
\left|\sigma_{n}\mathbb{P}(\frac{\eta_{n}-\mu_{n}}{\sigma_{n}}=x)-
\frac{1}{\sqrt{2\pi}}e^{-\frac{x^2}{2}}\right|=\frac{1}{2}\sqrt{n}\
\mathbb{P}(\frac{\eta_{n}-\frac{1}{4}n}{\frac{1}{2}\sqrt{n}}=x)-
\frac{1}{\sqrt{2\pi}}e^{-\frac{x^2}{2}}
$$
and for $S=\{\frac{\sqrt{n}}{2}|n=1,2\ldots\}$, we have $\left|\frac{1}{2}
\sqrt{n}\mathbb{P}(\eta_{n}=
\frac{1}{2}n)-\frac{1}{\sqrt{2\pi}}e^{-\frac{n}{4}}\right|\nrightarrow
0$, the key point being that $a_{n,k}$ flips between $-\binom{n}{k}$
and $3{n \choose k}$.

All results of this paper hold for $2$-noncrossing RNA structures,
i.e.~RNA secondary structures. This is a consequence of an analogous
analysis of their respective bivariate generating function. In this
case, however, no singular expansion is necessary as the generating
function itself can be used. They also give rise to put the asymptotic
results on RNA secondary structures of \cite{Schuster:98a} on a new level.
We can pass from computing exponential growth rates to computing
distributions for RNA secondary structures with specific properties.
To be precise we have for RNA secondary
structures

{\bf Theorem.} {\it Let ${\sf S}_{2}'(n,h)$ denote the number of RNA
secondary structures with exactly $h$ arcs. Then the random variable
$Y_{n}$ having distribution $\mathbb{P}(Y_{n}=h)={\sf
S}_{2}'(n,h)/{\sf S}_{2}(n)$ satisfies a central and local limit
theorem with mean $0.27639\,n$ and variance $0.04472\, n$. }

In particular the theorem predicts a sharp concentration of the
number of RNA secondary structures with $55.278\%$ unpaired bases
which agrees with the statistics of RNA secondary structures obtained
by folding algorithms
\cite{Zuker:79b,Schuster:98a,Waterman:86,Bauer:96,Tacker:94a,McCaskill:90a}.
Let us finally remark that much more holds: due to the determinant
formula for $k$-noncrossing matchings and the functional identity of
Lemma~\ref{L:func}, Section~\ref{S:func} our results can be generalized to
$k$-noncrossing RNA structures, where $k$ is arbitrary. Why this is of
interest can be seen in Fig.~\ref{F:4}. For higher $k$ the mean of the
central limit theorems for $k$-noncrossing RNA structures will shift towards
the maximum combinatorially possible number of arcs. We speculate that
each increase in $k$ will basically cut the distance to the maximum arc
number in half. This is work in progress.

The paper is structured as follows: In Section~\ref{S:rna} we
provide some background on $k$-noncrossing RNA structures and all
generating functions involved. In Section~\ref{S:func} we give a
functional equation for the bivariate generating function of ${\sf
S}_{3}'(n,h)$ via $3$-noncrossing matchings proved in
\cite{Reidys:07rna2}. We have included its proof in the appendix in
order to keep the paper self-contained. This functional identity
plays a key role in proving the central limit and local limit
theorem in Section~\ref{S:CL} and Section~\ref{S:LL}, respectively.
The central limit theorem is proved by analyzing the singular
expansion of analytic function of power series $\sum_{n \ge
0}\sum_{h\le \frac{n}{2}}{\sf S}_{3}'(n,h)w^hz^n$ and using transfer
theorems \cite{Flajolet:05,Gao:92,Handbook:95} and to prove the
local limit theorem, we use a theorem of Hwang \cite{Hwang:98} and
build on our proof of the central limit theorem.


\section{RNA structures}\label{S:rna}

\begin{figure}[ht]
\centerline{ \epsfig{file=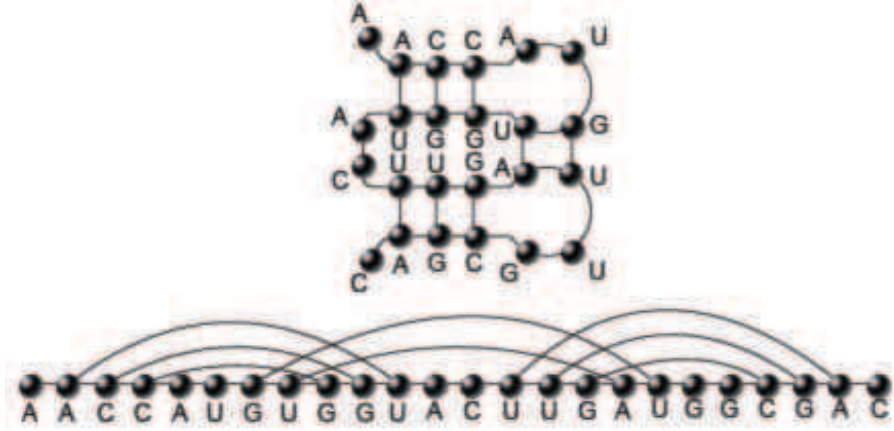,width=0.8\textwidth}\hskip15pt}
\caption{\small A $3$-noncrossing RNA structure, as a planar graph
(top) and as a diagram (bottom)} \label{F:5}
\end{figure}

Let us begin by illustrating the concept of RNA structures. Suppose we are
given the primary sequence
$$
{\bf A}{\bf A}{\bf C}{\bf C}{\bf A}{\bf U}{\bf G}{\bf U}{\bf G}{\bf G}
{\bf U}{\bf A}{\bf C}{\bf U}{\bf U}{\bf G}{\bf A}{\bf U}{\bf G}{\bf G}
{\bf C}{\bf G}{\bf A}{\bf C}  \ .
$$
Structures are combinatorial graphs over the labels of the nucleotides
of the primary sequence. These graphs can be represented in several ways.
In Figure~\ref{F:5} we represent a $3$-noncrossing RNA structure with loop-loop
interactions in two ways:
first we display the structure as a planar graph
and secondly as a diagram, where the bonds are drawn as arcs in the
positive half-plane.


In the following we will consider structures as diagram
representations of digraphs. A digraph $D_n$ is a pair of sets
$V_{D_n},E_{D_n}$, where $V_{D_n}= \{1,\dots,n\}$ and
$E_{D_n}\subset \{(i,j)\mid 1\le i< j\le n\}$. $V_{D_n}$ and
$E_{D_n}$ are called vertex and arc set, respectively. A
$k$-noncrossing digraph is a digraph in which all vertices have
degree $\le 1$ and which does not contain a $k$-set of arcs that are
mutually intersecting, i.e.
\begin{eqnarray}
\not\exists\,
(i_{r_1},j_{r_1}),(i_{r_2},j_{r_2}),\dots,(i_{r_k},j_{r_k});\quad & &
i_{r_1}<i_{r_2}<\dots<i_{r_k}<j_{r_1}<j_{r_2}<\dots<j_{r_k} \ .
\end{eqnarray}
We will represent digraphs as a diagrams (Figure~\ref{F:5}) by
representing the vertices as integers on a line and connecting any
two adjacent vertices by an arc in the upper-half plane. The direction
of the arcs is implicit in the linear ordering of the vertices and
accordingly omitted.

\begin{definition}\label{D:rna}
An RNA structure (of pseudo-knot type $k-2$), ${S}_{k,n}$, is a
digraph in which all vertices have degree $\le 1$, that does not
contain a $k$-set of mutually intersecting arcs and $1$-arcs,
i.e.~arcs of the form $(i,i+1)$, respectively. We denote the number
of RNA structures by ${\sf S}_k(n)$ and the number of RNA structures
with exactly $\ell$ isolated vertices and with $h$ arcs by ${\sf
S}_k(n,\ell)$ and ${\sf S}_k'(n,h)$, respectively. Note that ${\sf
S}_k'(n,h)={\sf S}_k(n,n-2h)$.
\end{definition}

Let $f_{k}(n,\ell)$ denote the number of $k$-noncrossing digraphs with
$\ell$ isolated points. We have shown in \cite{Reidys:07rna1} that
\begin{align}\label{E:ww0}
f_{k}(n,\ell)& ={n \choose \ell} f_{k}(n-\ell,0) \\
\label{E:ww1}
\det[I_{i-j}(2x)-I_{i+j}(2x)]|_{i,j=1}^{k-1} &=
\sum_{n\ge 1} f_{k}(n,0)\cdot\frac{x^{n}}{n!} \\
\label{E:ww2}
e^{x}\det[I_{i-j}(2x)-I_{i+j}(2x)]|_{i,j=1}^{k-1}
&=(\sum_{\ell \ge 0}\frac{x^{\ell}}{\ell!})(\sum_{n \ge
1}f_{k}(n,0)\frac{x^{n}}{n!})=\sum_{n\ge 1}
\left\{\sum_{\ell=0}^nf_{k}(n,\ell)\right\}\cdot\frac{x^{n}}{n!} \ .
\end{align}
In particular we obtain for
$k=2$ and $k=3$
\begin{equation}\label{E:2-3}
f_2(n,\ell)  =  \binom{n}{\ell}\,C_{(n-\ell)/2}\quad
\text{\rm and}\quad  f_{3}(n,\ell)=
{n \choose \ell}\left[C_{\frac{n-\ell}{2}+2}C_{\frac{n-\ell}{2}}-
      C_{\frac{n-\ell}{2}+1}^{2}\right] \ ,
\end{equation}
where $C_m$ denotes the $m$-th Catalan number.
The derivation of the generating function of $k$-noncrossing RNA structures,
given in Theorem~\ref{T:cool1} below uses advanced methods and novel
constructions of enumerative combinatorics due to Chen~{\it et.al.}
\cite{Chen:07a,Gessel:92a} and Stanley's mapping between matchings and
oscillating tableaux i.e.~families of Young diagrams in which any two
consecutive shapes differ by exactly one square.
The enumeration is obtained using the
reflection principle due to Gessel and Zeilberger \cite{Gessel:92a} and
Lindstr\"om \cite{Lindstroem:73a} combined with an inclusion-exclusion
argument in order to eliminate the arcs of length $1$. In
\cite{Reidys:07rna1} generalizations to restricted (i.e. where arcs of
the form $(i,i+2)$ are excluded) and circular RNA structures are given.
The following theorem provides all data on numbers of $k$-noncrossing
RNA structures with $h$ arcs and the numbers of all $k$-noncrossing RNA
structures.
\begin{theorem}\cite{Reidys:07rna1}\label{T:cool1}
Let $k\in\mathbb{N}$, $k\ge 2$, let $C_m$ denote the $m$-th Catalan number
and $f_k(n,\ell)$ be the number of $k$-noncrossing digraphs over $n$
vertices with exactly $\ell$ isolated vertices. Then the
number of RNA structures with $\ell$ isolated vertices,
${\sf S}_k(n,\ell)$, is given by
\begin{equation}\label{E:da}
{\sf S}_k(n,\ell) = \sum_{b=0}^{(n-\ell)/2}
                   (-1)^b\binom{n-b}{b}f_k(n-2b,\ell)  \  ,
\end{equation}
where $f_k(n-2b,\ell)$ is given by the generating function in
eq.~{\rm (\ref{E:ww1})}.
Furthermore the number of $k$-noncrossing RNA structures, ${\sf S}_k(n)$
is
\begin{equation}\label{E:sum}
{\sf S}_k(n)
=\sum_{b=0}^{\lfloor n/2\rfloor}(-1)^{b}{n-b \choose b}
\left\{\sum_{\ell=0}^{n-2b}f_{k}(n-2b,\ell)\right\}
\end{equation}
where $\{\sum_{\ell=0}^{n-2b}f_{k}(n-2b,\ell)\}$ is given by the
generating function in eq.~{\rm (\ref{E:ww2})}.
\end{theorem}

In principle, Theorem~\ref{T:cool1} contains all information about
the numbers of $k$-noncrossing RNA structures. However, due to the
inclusion-exclusion structure of its coefficients it is however
difficult to interpret and to express their behavior for large $n$.
Subsequent asymptotic analysis \cite{Reidys:07rna2}
produced the following simple formula
\begin{theorem}\cite{Reidys:07rna2}\label{T:asy3}
The number of $3$-noncrossing RNA structures is asymptotically given by
\begin{eqnarray*}
\label{E:konk3}
{\sf S}_3(n) & \sim & \frac{10.4724\cdot 4!}{n(n-1)\dots(n-4)}\,
\left(\frac{5+\sqrt{21}}{2}\right)^n \ .\\
\end{eqnarray*}
\end{theorem}


\section{A functional equation}\label{S:func}


We have shown in the introduction that the bivariate generating
function is the key to prove the central and local limit theorems.
The following lemma, whose proof is given in the appendix, rewrites
this bivariate generating function as a composition of two
``simple'' functions. This is crucial for the singularity analysis
insofar as we can use a phenomenon known as persistence of the
singularity of the ``outer'' function (the {\it supercritical case})
\cite{Flajolet:05}. It basically means that the type of the singularity
is determined by the generating function of $k$-noncrossing matchings.
\begin{lemma}\label{L:func}\cite{Reidys:07rna2}
Let $x$ be an indeterminant over $\mathbb{R}$ and $w\in\mathbb{R}$ a
parameter. Let $\rho_k(w)$ denote the radius of convergence
of the power series $\sum_{n\ge 0} [\sum_{h\le n/2} {\sf S}'_k(n,h)
w^{2h}] x^n$. Then for $ \vert x \vert < \rho_k(w)$
\begin{equation}\label{E:rr}
\sum_{n\ge 0} \sum_{h\le n/2} {\sf S}_k'(n,h) w^{2h} x^n =
\frac{1}{w^2x^2-x+1} \sum_{n\ge 0} f_k(2n,0)
\left(\frac{wx}{w^2x^2-x+1}\right)^{2n}
\end{equation}
holds. In particular we have for $w=1$
\begin{equation}\label{E:oha}
\sum_{n\ge 0} {\sf S}_k(n) z^{n} =\frac{1}{z^2-z+1}\, \sum_{n\ge
0} f_k(2n,0) \left(\frac{z}{z^2-z+1}\right)^{2n}
\end{equation}
for $z\in \mathbb{C}$ with $\vert z\vert <\rho_{k}(1)$.
\end{lemma}
To keep the paper selfcontained we give the proof of Lemma~\ref{L:func}
in the Appendix.
While ~(\ref{E:rr}) can only be proved on the level of formal
power-series for real variables, complex analysis i.e.~the
interpretation of these generating functions as analytic functions
allows to extend the equality to arbitrary complex variables.

\begin{lemma}\label{L:bivariate}
Suppose $\epsilon>0$, $k\in\mathbb{N}$, $k\ge 2$ and
$w=e^{\frac{s}{2}}$, where $\vert s\vert <\epsilon$ and
$\varphi_{n,k}(s)=\sum_{h\le n/2}{\sf S}_k'(n,h) e^{hs}$. Let
$\rho_k(s)\in \mathbb{R^+}$ denote the radius of convergence of
$\sum_{n\ge 0} \varphi_{n,k}(s)z^n$ parameterized by $s$. Then we
have
\begin{equation}\label{E:kl}
\forall\, s,z\in\mathbb{C};\
\vert s\vert <\epsilon, \vert z\vert< \rho_k(s);\quad
\sum_{n\ge 0}\varphi_{n,k}(s)z^n= \frac{1}{e^{s}z^2-z+1} \sum_{n\ge
0} f_k(2n,0) \left(\frac{e^{\frac{s}{2}}z}{e^{s}z^2-z+1}\right)^{2n}
\ .
\end{equation}
Furthermore $\sum_{n\ge 0}\varphi_{n,3}(s)z^n$ has an analytic
continuation, $\Xi_3(z,s)$. For $\epsilon$ sufficiently small and
$\vert s\vert <\epsilon$, $\Xi_3(z,s)$ has exactly $6$
singularities, $4$ of which have distinct moduli.
\end{lemma}
\begin{proof}
We first prove eq.~(\ref{E:kl}). For this purpose we observe that
\begin{equation}\label{E:G}
\forall\, \vert s\vert<\epsilon,\; \vert z\vert < \rho_k(s) \qquad
G(z,s)=\frac{1}{e^{s}z^2-z+1} \sum_{n\ge 0} f_k(2n,0)
\left(\frac{e^{\frac{s}{2}}z}{e^{s}z^2-z+1}\right)^{2n}
\end{equation}
considered as a power series in $e^{\frac{1}{2}s}$ is analytic in a
neighborhood of $s=0$, since $G(z,0)$ is analytic for $\vert z \vert
<\rho_{k}(0)$. In addition, we can interpret $\sum_{n\ge
0}\varphi_{n,k}(s)z^n$ as a power series in $e^{\frac{1}{2}s}$:
\begin{equation}\label{E:ana-s}
\sum_{n\ge 0}\left(\sum_{h\le n/2}{\sf S}_k'(n,h)e^{hs}\right)z^n
=\sum_{h\ge 0}\left(\sum_{n\ge 2h}{\sf S}_k'(n,h)z^n\right)
\left(e^{s}\right)^h =\sum_{h\ge 0}\psi_{h}(z)
\left(e^{\frac{1}{2}s}\right)^{2h} \ .
\end{equation}
Therefore $G(z,s)$ and the power series $\sum_{n\ge 0}\varphi_{n,k}(s)z^n$ are
analytic in the indeterminant $e^{\frac{1}{2}s}$ in an $\epsilon$-disc
centered at $0$. Lemma~\ref{L:func} implies that for
$s\in ]-\epsilon,\epsilon[$ the analytic functions $G(z,s)$ and
$\sum_{n\ge 0}\varphi_{n,k}(s)z^n$ are equal. Since any two functions that are
analytic at $0$ and that coincide on the interval $]-\epsilon,\epsilon[$
are identical, we obtain
\begin{equation}\label{E:gleich}
\forall\, \vert s\vert<\epsilon,\; \vert z\vert < \rho_k(s)\qquad
G(z,s)=\sum_{n\ge 0}\varphi_{n,k}(s)z^n \ .
\end{equation}
{\it Claim $1$.} Suppose $\vert s\vert <\epsilon$. Then $\sum_{n\ge 0}
\varphi_{n,3}(s)z^n$ has an analytic continuation, $\Xi_3(z,s)$, which has
exactly $6$ singularities $4$ of which have distinct moduli.
\\
In order to prove Claim $1$ we observe that the power series
$\sum_{n\ge 0} f_3(2n,0) y^{n}$ has the analytic continuation $\Psi(y)$
(obtained by MAPLE sumtools) given by
\begin{equation}\label{E:psi}
\Psi(y)= \frac{-(1-16y)^{\frac{3}{2}}
P_{\frac{3}{2}}^{-1}(-\frac{16y+1}{16y-1})} {16\, {y}^{\frac{5}{2}}}
\ ,
\end{equation}
where $P_{\nu}^{m}(x)$ denotes the Legendre Polynomial of the first kind
with the parameters $\nu=\frac{3}{2}$ and $m=-1$.
According to eq.~(\ref{E:gleich}) we have
\begin{equation}\label{E:G2}
\sum_{n\ge 0}\varphi_{n,k}(s)z^n=\frac{1}{e^{s}z^2-z+1}
\sum_{n\ge 0} f_k(2n,0) \left(\frac{e^{\frac{s}{2}}z}{e^{s}z^2-z+1}\right)^{2n}
\end{equation}
which implies that $\sum_{n\ge 0}\varphi_{n,3}(s)z^n$ has the analytic
continuation
\begin{equation}
\forall\, \vert s\vert<\epsilon,\;\qquad \Xi_3(z,s)=
\frac{1}{e^sz^2-z+1}\,
\Psi\left(\frac{e^{\frac{1}{2}s}z}{e^sz^2-z+1}\right)^2\ .
\end{equation}
In particular for $s=0$, $\Xi_3(z,0)$ is the analytic continuation
of the power series $\sum_{n\ge 0} {\sf S}_3(n) z^n$. We proceed by
showing that $\Xi_3(z,s)$ has exactly $6$ singularities and 4 of
them have different moduli in $\mathbb{C}$ parameterized by $s$. Two
singularities are given by the roots of $e^sz^2-z+1$ are
$\zeta_1(s)=\frac{1-\sqrt{1-4e^s}}{2e^s}$ and
$\zeta_2(s)=\frac{1+\sqrt{1-4e^s}}{2e^s}$. Observe that $\vert
\zeta_1(0)\vert=\vert\zeta_2(0) \vert=1$ and polynomial $e^sz^2-z+1$
depends continuously on $e^{\frac{s}{2}}$, therefore $\zeta_1(s)$
and $\zeta_2(s)$ could potentially have equal modulus for $\vert
s\vert<\epsilon$. The remaining $4$ singularities are induced by the
the unique dominant singularity $\alpha_1=\frac{1}{16}$ of analytic
function $\Psi(y)$. The function $\Psi(y)$ has three
singularities, two of them $\alpha_1=\frac{1}{16}$ and
$\alpha_2=+\infty$ are branch points and the other $\alpha_3=0$ is a
removable singularity. The function
$g(z)=\left(\frac{e^{\frac{s}{2}}z}{e^sz^2-z+1}\right)^2$ with
$g(0)=0$ has a radius of convergence of $1$ as $s$ tends to $0$.
Therefore the singularity type only depends on
$\Psi(y)$ (this is the {\it supercritical case} in \cite{Flajolet:05}).
The singularity $\alpha_1=\frac{1}{16}$ gives rise to the equations
$$
0=e^sz^2-(1+4e^{\frac{1}{2}s})z +1 \quad \text{\rm and} \quad
0=e^sz^2+(4e^{\frac{1}{2}s}-1)z+1
$$
and setting $\mu_+(s)=1+ 4e^{\frac{1}{2}s}$, $\mu_-(s)=1- 4e^{\frac{1}{2}s}$
and $\theta(s)=\sqrt{12e^s+8e^{\frac{1}{2}s}+1}$ its roots
are given by
$$
\zeta_3(s)=\frac{\mu_+(s)-\theta(s)}{2e^s}, \
\zeta_4(s)=\frac{\mu_+(s)+\theta(s)}{2e^s},\
\zeta_5(s)=\frac{\mu_-(s)+\theta(s)}{2e^s}\ \text{\rm and}\
\zeta_6(s)=\frac{\mu_-(s)-\theta(s)}{2e^s},
$$
respectively. Observe that for $\vert s \vert <\epsilon$,
$e^{\frac{s}{2}}$ is in a neighborhood of $1$ over $\mathbb{C}$,
hence $\theta(s)\ne 0$. That leads to $4$ distinct roots
$\zeta_{3}(s),\zeta_{4}(s),\zeta_{5}(s),\zeta_{6}(s)$ over $\vert
s\vert <\epsilon$, all of them have distinct moduli for $s$ being a
sufficiently small neighborhood of 0. Indeed, for $s=0$ we have $4$
distinct real valued roots
$$
\zeta_3(0)=\frac{5-\sqrt{21}}{2}, \
\zeta_4(0)=\frac{5+\sqrt{21}}{2},\ \zeta_5(0)=\frac{-3+\sqrt{5}}{2},
\ \text{\rm and}\ \zeta_6(0)=\frac{-3-\sqrt{5}}{2}
$$
and the polynomials $e^sz^2-(1+4e^{\frac{1}{2}s})z
+1$, $e^sz^2+(4e^{\frac{1}{2}s}-1)z+1$ and $e^sz^2-z+1$ depend
continuously on the parameter $e^{\frac{1}{2}s}$, whence Claim $1$
and the lemma follows.
\end{proof}

\section{The central limit theorem}\label{S:CL}

In this section we prove a central limit theorem for the numbers of
$3$-noncrossing RNA structures with $h$ arcs.
We will analyze for fixed but arbitrary $n$ the
distribution of ${\sf S}_3'(n,h)$. Let us first prepare some
methods and results used in the proof of Theorem~\ref{T:gauss}.
$[z^n]\,f(z)$ denotes the coefficient of $z^n$ in the power series
expansion of $f(z)$ around $0$.
The scaling property of Taylor coefficients
\begin{equation}\label{E:scaling}
\forall \,\gamma\in\mathbb{C}\setminus 0;\quad
[z^n]f(z)=\gamma^n [z^n]f(\frac{z}{\gamma}) \ ,
\end{equation}
shows that w.l.o.g.~any singularity analysis can be reduced to the
case where $1$ is the dominant singularity. We will be interested
in the behavior of an analytic function ``locally'', i.e.~around a
certain singularity $\rho$. For this purpose we use the notation
\begin{equation}\label{E:genau}
f(z)=O\left(g(z)\right) \ \text{\rm as $z\rightarrow \rho$}\quad
\Longleftrightarrow \quad f(z)/g(z) \ \text{\rm is bounded as
$z\rightarrow \rho$}
\end{equation}
and if we write $f(z)=O(g(z))$ it is implicitly assumed that $z$
tends to a (unique) singularity.
Given two numbers $\phi,R$, where $R>\vert \rho \vert>0$ and
$0<\phi<\frac{\pi}{2}$ and $\rho \in\mathbb{C}$ the open domain
$\Delta_{\rho}(\phi,R)$ is defined as
\begin{equation}
\Delta_{\rho}(\phi,R)=\{ z\mid \vert z\vert < R, z\neq \rho,\, \vert
{\rm Arg}(z-\rho)\vert
>\phi\}
\end{equation}
A domain is a $\Delta_{\rho}$-domain if it is of the form
$\Delta_{\rho}(\phi,R)$ for some $R$ and $\phi$. A function is
$\Delta_{\rho}$-analytic if it is analytic in some
$\Delta_{\rho}$-domain. We use $U(a,r)=\{z\in \mathbb{C} \mid
\vert z-a\vert<r\}$ to denote the open neighborhood of $a$ in
$\mathbb{C}$. Via the following theorem we can
extract the coefficients of analytic functions provided these functions
satisfy certain ``local'' properties.
\begin{theorem}\label{T:transfer1}\cite{Flajolet:05}
Let $r\in\mathbb{Z}_{\ge 0}$ and $f(z,e^s)$ be a $\Delta_{\rho(s)}$-analytic
function parameterized by $s$, which satisfies in the intersection of a
neighborhood of $\rho(s)$ with its $\Delta_{\rho(s)}$-domain
\begin{equation}\label{T:transfer2}
f(z,e^s) = b_0(s)+b_1(s)(z-\rho(s))+A(s) \
(\rho(s)-z)^{r}\ln^{}\left(\frac{1}{\rho(s)-z}\right) + R(z,s)
\end{equation}
where $A(s),b_0(s),b_1(s)$ are analytic in $\vert s\vert<
\epsilon$ and $\vert R(z,s)\vert \le c \, \vert \rho(s)-z \vert$
for some absolute constant $c\in \mathbb{C}$. That is we have
$f(z,e^s)=O((\rho(s)-z)^{r}\ln(\frac{1}{\rho(s)-z}))$ with uniform
error bound as $s$ in a neighborhood of $0$. Then we have
\begin{equation}
[z^n]f(z,e^s)= A(s)\,  (-1)^r\frac{r!}{n(n-1)\dots(n-r)}
\left(1-O(\frac{1}{n})\right)\quad \text{\it for some
$A(s)\in\mathbb{C}$}\ ,
\end{equation}
where the error term is again uniform for $s$ from a neighborhood
of origin, i.e.~$R(s)\le c\, \vert s\vert$, where $c>0$.
\end{theorem}
{\bf Remark.}
The equivalence between eq.~(\ref{T:transfer2}) and $f(z,e^s)=
O((\rho(s)-z)^{r}\ln(\frac{1}{\rho(s)-z}))$ for $r\in\mathbb{Z}_{\ge 0}$
can be seen as follows: by definition of $f(z,e^s)=
O((\rho(s)-z)^{r}\ln(\frac{1}{\rho(s)-z}))$ there exist $A(z,s)$ and
$B(z,s)$, such that
$f(z,e^s)=B(z,s)+A(z,s)(\rho(s)-z)^{r}\ln(\frac{1}{\rho(s)-z})$,
where $A(z,s)$ and $B(z,s)$ are analytic in a neighborhood of
$\rho(s)$. Taylor expansion of $A(z,s)$ and $B(z,s)$ at $z=\rho(s)$
produces
\begin{align*}
f(z,s)&=B(z,s)+A(z,s)(\rho(s)-z)^{r}\ln\left(\frac{1}{\rho(s)-z}\right)\\
&=b_{0}(s)+b_{1}(s)(z-\rho(s))+\cdots+\left(a_0(s)+a_{1}(z-\rho(s))+
\cdots\right)
(\rho(s)-z)^{r}\ln\left(\frac{1}{\rho(s)-z}\right)\\
&=b_{0}(s)+b_{1}(s)(z-\rho(s))+a_{0}(s)(\rho(s)-z)^{r}
\ln\left(\frac{1}{\rho(s)-z}\right)+R(z,s)
\end{align*}
where
$R(z,s)=O((\rho(s)-z)^{r+1}\ln\left(\frac{1}{\rho(s)-z}\right))$.
For $r \in \mathbb{Z}_{\ge 0}$,
$\frac{|R(z,s)|}{|\rho(s)-z|}=O(|\rho(s)-z|^{r}
\ln\left(\frac{1}{|\rho(s)-z|}\right))$
is bounded by an absolute constant as $z$ tends to $\rho(s)$. That
implies the error bound is uniform.

The next Theorem is a classic result on limit distributions which
allows us to prove our main result via characteristic functions
i.e.~explicitly by showing $\lim_{n \rightarrow
\infty}\varphi_{n}(t)=\varphi(t)$ for any $t \in (-\infty,\infty)$.
\begin{theorem}\label{T:K}{\bf (L\'{e}vy-Cram\'{e}r)}
Let $\{\xi_{n}\}$ be a sequence of random variables and let
$\{\varphi_{n}(x)\}$ and $\{F_{n}(x)\}$ be the corresponding
sequences of characteristic and distribution functions. If there
exists a function $\varphi(t)$, such that $\lim_{n \rightarrow
\infty}\varphi_{n}(t)=\varphi(t)$ uniformly over an arbitrary
finite interval enclosing the origin, then there exists a random
variable $\xi$ with distribution function $F(x)$ such that
$$F_{n}(t)\Longrightarrow F(x)$$
uniformly over any finite or infinite interval of continuity of
$F(x)$.
\end{theorem}


We now consider the random variable $X_n$ having the distribution
$\mathbb{P}(X_n=h)={\sf S}_3'(n,h)/{\sf S}_3(n)$, where
$h=0,1,\ldots \lfloor \frac{n}{2}\rfloor$. The key point in the
proof of Theorem~\ref{T:gauss} is to compute the coefficients of
the bivariate generating function whose variable, $s$ is
considered as a parameter. Intuitively the particular distribution
is a result of how the singularity shifts as a function of this
parameter. As a result the proof is somewhat ``non-probabilistic''
and has two distinct parts: {\sf (a)} the analytic combinatorics
of the bivariate generating function and {\sf (b)} the computation
of the characteristic function with subsequent application of the
{L\'{e}vy-Cram\'{e}r} Theorem.
\begin{theorem}\label{T:gauss}
The random variable $\frac{X_{n}-\mu n}{\sqrt{\sigma^2 n}}$ has asymptotically
normal distribution with parameter $(0,1)$, i.e.
\begin{equation}\label{E:converge}
\lim_{n\to\infty}\mathbb{P}\left(\frac{X_n-\mu n}{\sqrt{\sigma^2 n}}< x
\right)  =  \frac{1}{\sqrt{2\pi}}\int_{-\infty}^{x}
                                          e^{-\frac{1}{2}t^2} dt
\end{equation}
and $\mu,\sigma^2$ are given by
\begin{equation}\label{E:concrete0}
\mu  =  -\frac{-\frac{3}{2}+\frac{13}{42}\sqrt{21}}{\frac{5}{2}-
                   \frac{1}{2}\sqrt{21}}=0.39089  \quad \text{\rm and}
                   \quad
\sigma^2  =  \mu^2-\frac{1-\frac{94}{441}\sqrt{21}}
{\frac{5-\sqrt{21}}{2}}=0.041565  \ .
\end{equation}
\end{theorem}
\begin{proof}
We set $w=e^{\frac{1}{2}s}$ and $\varphi_{n,3}(s)=\sum_{h\le
n/2}{\sf S}_3' (n,h)e^{hs}$. Since
\begin{equation}\label{E:phi}
\sum_{n\ge 0}\varphi_{n,3}(s)z^n = \sum_{n\ge 0}\left(\sum_{h\le
n/2}{\sf S}_3'(n,h)e^{hs}\right)z^n \ ,
\end{equation}
we can consider the double generating function $\sum_{n\ge 0}\sum_{h\le n/2}
{\sf S}_3'(n,h)w^{2h}z^n$ as a power series in the complex indeterminant $z$,
parameterized by $s$. \\
{\it Claim $1$.}
\begin{equation}
\Psi(z)=O\left((1-16z)^4\ln\left(\frac{1}{1-16z}\right)\right)\quad
\ \text{\rm holds uniformly for }\forall\,
z\in\Delta_{\frac{1}{16}}(\phi,R)\cap U(\frac{1}{16},\epsilon);
\end{equation}
$\Psi(z)$ is $\Delta_{\frac{1}{16}}(\phi,R)$-analytic and has the
singular expansion $(1-16z)^4\ln\left(\frac{1}{1-16z}\right)$ in the
intersection of $U(\frac{1}{16},\epsilon)$ with the
$\Delta_{\frac{1}{16}}-$domain, where $\Delta_{r}(\phi,R)=\{z\big|
|z|<R, z\ne r,\vert{\rm Arg}(z-r)\vert>\phi \}$ for some $R>r$.
First $\Delta_{\frac{1}{16}}(\phi,R)$-analyticity of the function
$(1-16z)^4\ln\left(\frac{1}{1-16z}\right)$ is obvious. We proceed by
proving that $(1-16z)^4\ln\left(\frac{1}{1-16z}\right)$ is the
singular expansion of $\Psi(z)$. The above mentioned scaling
property of Taylor coefficients allows us to consider the power
series $\sum_{n \ge 0}f_{3}(2n,0)(\frac{z} {16})^n$ over the
$\Delta$-domain $\Delta_{1}(\phi,R)$ for some $R>1$. Using the
notation of falling factorials $(n-1)_4=(n-1)(n-2)(n-3)(n-4)$ we
observe
$$
f_3(2n,0)=C_{n+2}C_{n}-C_{n+1}^2= \frac{1}{(n-1)_4}
\frac{12(n-1)_4(2n+1)}{(n+3)(n+1)^2(n+2)^2}\,
\binom{2n}{n}^2 \ .
$$
With this expression for $f_3(2n,0)$ we arrive at the formal identity
\begin{eqnarray*}
\sum_{n\ge 5}16^{-n}f_3(2n,0)z^n  & = &
O(\sum_{n\ge 5}
\left[16^{-n}\,\frac{1}{(n-1)_4}
\frac{12(n-1)_4(2n+1)}{(n+3)(n+1)^2(n+2)^2}\,
\binom{2n}{n}^2-\frac{4!}{(n-1)_4}\frac{1}{\pi}\frac{1}{n}\right]z^n \\
& & + \sum_{n\ge 5}\frac{4!}{(n-1)_4}\frac{1}{\pi}\frac{1}{n}z^n) \ ,
\end{eqnarray*}
where $f(z)=O(g(z))$ denotes that the limit $f(z)/g(z)$ is bounded
for $z\rightarrow 1$, eq.~(\ref{E:genau}). It is clear that the
error bound below
\begin{eqnarray*}
& &\sum_{n\ge 5}\left[16^{-n}\,\frac{1}{(n-1)_4}
\frac{12(n-1)_4(2n+1)}{(n+3)(n+1)^2(n+2)^2}\,
\binom{2n}{n}^2-\frac{4!}{(n-1)_4}\frac{1}{\pi}\frac{1}{n}\right]z^n  \\
&\sim & \sum_{n\ge 5} \left[16^{-n}\,\frac{1}{(n-1)_4}
\frac{12(n-1)_4(2n+1)}{(n+3)(n+1)^2(n+2)^2}\,
\binom{2n}{n}^2-\frac{4!}{(n-1)_4}\frac{1}{\pi}\frac{1}{n}\right]
 <\kappa
\end{eqnarray*}
holds uniformly for $z$ in $\Delta_{1}(\phi,R)\cap U(1,\epsilon)$
and some absolute $\kappa< 0.0784$. Therefore we can conclude
\begin{equation}
\sum_{n\ge 5}16^{-n}f_3(2n,0)z^n=
O(\sum_{n\ge 5}\frac{4!}{(n-1)_4}\frac{1}{\pi}\frac{1}{n}z^n) \ .
\end{equation}
We proceed by interpreting the power series on the rhs, observing
\begin{equation}
\forall\, n\ge 5\, ; \qquad
[z^n]\left((1-z)^4\,\ln\frac{1}{1-z}\right)=
\frac{4!}{(n-1)\dots (n-4)}\frac{1}{n} \, ,
\end{equation}
whence $\left((1-z)^4\,\ln\frac{1}{1-z}\right)$ is
the unique analytic continuation of $\sum_{n\ge 5}\frac{4!}{(n-1)_4}
\frac{1}{\pi}\frac{1}{n}z^n$.
Using the scaling property of Taylor coefficients
$[z^n]f(z)=\gamma^n [z^n]f(\frac{z}{\gamma})$
\begin{equation}\label{E:isses}
\Psi(z) =O\left((1-16z)^4\ln\left(\frac{1}{1-16z}\right)\right)
\mbox{ holds uniformly for }\forall\,
z\in\Delta_{\frac{1}{16}}(\phi,R)\cap U(\frac{1}{16},\epsilon)
\end{equation}
Therefore we have proved that $(1-16z)^{4}\ln(\frac{1}{1-16z})$ is
the singular expansion of $\Psi(z)$ at $z=\frac{1}{16}$, whence
Claim $1$. Our next step consists in verifying that when passing from
$\Psi(z)$ to the bivariate generating function
$\Psi(z,s)= \Psi((\frac{wz}{w^2z^2-z+1})^2)$,
then there exists a singular expansion of the form
${O}\left((1-\frac{z}{\rho_3(s)})^4\ln(\frac{1}{1-\frac{z}
{\rho_3(s)}})\right)$, parameterized in $s$.\\
{\it Claim $2$}. Let $0<\epsilon<1$, then for any $\vert
s\vert<\epsilon$ and $z\in \Delta_{\rho_3(s)}(\phi,R)$, we have
$\Psi(z,s)
={O}\left((1-\frac{z}{\rho_3(s)})^4\ln(\frac{1}{1-\frac{z}
{\rho_3(s)}})\right)$, and the error bound is uniform for $s$ in a
neighborhood of $0$. \\
To prove the claim we first observe that Claim $1$ implies
\begin{equation}
\Psi(z)=\kappa(1-16z)^4\ln\left(\frac{1}{1-16z}\right)+R(z)
\end{equation}
for some absolute constant $\kappa$ and $R(z)$ is the uniform error
bound for $z\in \Delta_{\frac{1}{16}}(\phi,R)\cap
U(\frac{1}{16},\epsilon)$. I.e.~For $z\in
\Delta_{\frac{1}{16}}(\phi,R)\cap U(\frac{1}{16},\epsilon)$, there
exists some absolute constant $c$, such that $|R(z)|\le c\cdot
|1-16z|$ holds. According to Lemma~\ref{L:bivariate} we have
\begin{align*}
\Xi_3(z,s) &=\frac{1}{e^sz^2-z+1}\ O
\,\left(\left(1-16(\frac{e^{\frac{1}{2}s}z}
{e^sz^2-z+1})^2\right)^4\ln\frac{1}{\left(1-16(
\frac{e^{\frac{1}{2}s}z}{e^sz^2-z+1})^2\right)}\right) \\
&=\frac{\kappa}{e^{s}z^2-z+1}
\left(1-16(\frac{e^{\frac{1}{2}s}z}{e^{s}z^2-z+1})^2\right)^4
\ln\left(\frac{1}{1-16(\frac{e^{\frac{1}{2}s}z}{e^{s}z^2-z+1})^2}\right)+
R(z,e^s) \ .
\end{align*}
We expand
$\left(1-16(\frac{e^{\frac{1}{2}s}z}{e^{s}z^2-z+1})^2\right)^4
\ln\left(\frac{1}{1-16(\frac{e^{\frac{1}{2}s}z}{e^{s}z^2-z+1})^2}\right)$
around $z=\rho_{3}(s)$, where $\rho_3(s)$ is the solution of
$\frac{ze^{\frac{1}{2}s}}{e^sz^2-z+1} =\frac{1}{4}$ of minimal
modulus. Lemma~\ref{L:bivariate} implies that
$$
\rho_3(s)=\frac{4e^{\frac{1}{2}s}+1-\sqrt{12e^s+8e^{\frac{1}{2}s}+1}}{2e^s}
$$
is the unique dominant singularity. As a function in $s$ we have
$\rho_3'(0)=-\frac{3}{2}+\frac{13}{42}\sqrt{21}\neq 0$. The term
$\sqrt{12e^s+8e^{\frac{1}{2}s}+1}$ in $\rho_{3}(s)$ produces two
branching points parameterized by $s$. i.e.
$w=e^{\frac{1}{2}s}=-\frac{1}{6}$ and
$w=e^{\frac{1}{2}s}=-\frac{1}{2}$, or equivalently $s=2\ln\frac{1}{2}+2\pi
i$ and $s=2\ln\frac{1}{6}+2\pi i$, respectively. The interval
between $2\ln\frac{1}{6}+2\pi i$ and $2\ln\frac{1}{2}+2\pi i$
divides the complex plane of $s$ into two analytic branches. For any
$0 <\epsilon < \min\{\vert 2\ln\frac{1}{2}+2\pi i\vert, \vert
2\ln\frac{1}{6}+2\pi i\vert\}=6.4343$, the region $\vert s \vert <
\epsilon$ is disjoint to the interval
$[(2\ln\frac{1}{6},2\pi),(2\ln\frac{1}{2},2\pi)]$. Therefore
$\rho_{3}(s)$ is analytic for $\vert s\vert<\epsilon$. We next
consider $q(z,s)=1-16(\frac{e^{\frac{1}{2}s}z} {e^sz^2-z+1})^2$ as a
function of $z$ and compute the Taylor expansion at $\rho_3(s)$.
$$
q(z,s)=\alpha(\rho_3(s)-z)+{O}(z-\rho_3(s))^2
$$
and setting $\alpha=\frac{\sqrt{21}}{5-\sqrt{21}}$
\begin{align*}
\frac{1}{e^sz^2-z+1}\, \left[q(z,s)^4\ln\frac{1} {q(z,s)}\right]
&= \frac{(\alpha(\rho_3(s)-z)+{O}(z-\rho_3(s))^2)^4\ln\frac{1}
{\alpha(\rho_3(s)-z)+{O}(z-\rho_3(s))^2}}{e^s(z-\rho_3(s))^2+(2\rho_3(s)e^s-1)
(z-\rho_3(s))-3\rho_3(s)^2e^s+\rho_3(s)+1}\\
&= \frac{\left([\alpha+O(z-\rho_3(s))](\rho_3(s)-z)^4
\ln\frac{1}{[\alpha+O(z-\rho_3(s))](\rho_3(s)-z)}
\right)}{O(z-\rho_3(s))-3\rho_3(s)^2+\rho_3(s)+1} \\
&=O(\rho_{3}(s)-z)^4\ln\left(\frac{1}{\rho_{3}(s)-z}\right) \ .
\end{align*}
According to Theorem~\ref{T:transfer1} for $r=4$, we obtain the
error term in the expansion of $\frac{1}{e^sz^2-z+1}\,
\left[q(z,s)^4\ln\frac{1} {q(z,s)}\right]$ is uniform for $s$ in a
neighborhood of $0$. We observe that the resulting error bound for
$\Xi_{3}(z,s)$ is the sum $R(z,e^s)+R_{1}(z,e^s)$, where
$$\vert R(z,e^s)\vert \le
c\cdot\big|1-16\left(\frac{e^{\frac{1}{2}s}z}{e^{s}z^2-z+1}\right)\big|
=O(\rho_{3}(s)-z) \ .
$$
Therefore the error bound for the expansion of bivariate $\Xi_{3}(z,s)$
is uniform and Claim $2$ is proved. We proceed by using the scaling property
of Taylor coefficients $[z^n]f(z)= \gamma^n [z^n]f(\frac{z}{\gamma})$ and apply
Theorem~\ref{T:transfer1}.
Via Theorem~\ref{T:transfer1} we obtain the key information about
the coefficients of $\Xi_3(z,s)$ which allows us to substitute
$\varphi_{n,3}(\frac{it}{\sigma_n})$ in eq.~(\ref{E:well}) below:
\begin{equation}\label{E:4}
[z^n]\,\Xi_3(z,s)= K(s) \,\frac{4!}{n(n-1)\dots(n-4)}\,
\left(\rho_3(s)^{-1}\right)^n \left(1-O(\frac{1}{n})\right)\quad
\text{\rm for some $K(s)\in\mathbb{C}$}\ ,
\end{equation}
where the error term is again uniform for $s$ from a neighborhood of
origin.  \\
Suppose we are given the random variable (r.v.) $\xi_n$ with mean $\mu_n$
and variance $\sigma_n^2$. We consider the rescaled r.v.~$\eta_{n}=(\xi_{n}
- \mu_{n})\sigma^{-1}_{n}$ and the characteristic function of $\eta_{n}$:
\begin{equation}
f_{\eta_n}(t)=\mathbb{E}[e^{it\eta_{n}}]=
 \mathbb{E}[e^{it\frac{\xi_{n}}{\sigma_{n}}}]
                                     e^{-i\frac{\mu_{n}}{\sigma_{n}}t} \ .
\end{equation}
In particular, for $\xi_n=X_n$ we obtain, substituting the term
$\mathbb{E}[e^{it\eta_{n}}]$
\begin{equation}
f_{X_n}(t)=\left(
\sum_{h=0}^{n}\frac{{\sf S}_3'(n,h)}{{\sf S}_3(n)}
e^{it\frac{h}{\sigma_{n}}}\right)\,
e^{-i\frac{\mu_{n}}{\sigma_{n}}t}
 \ .
\end{equation}
Since $\varphi_{n,3}(s)=\sum_{h\le n/2}{\sf S}_3'(n,h)e^{hs}$, we can interpret
${\sf S}_3(n)=\sum_{h\le n/2}{\sf S}_3'(n,h)$ as $\varphi_{n,3}(0)$ and
$\varphi_{n,3}(\frac{it}{\sigma_n})=
\sum_{h\le n/2}{\sf S}_3'(n,h)e^{h\frac{it}{\sigma_n}}$, respectively.
Therefore we have
\begin{equation}\label{E:well}
f_{X_n}(t)=\frac{1}{\varphi_{n}(0)}\,\varphi_{n}(\frac{it}{\sigma_{n}})\,
e^{-i\frac{\mu_{n}}{\sigma_{n}}t} \ .
\end{equation}
For $\vert s\vert < \epsilon$, eq.~(\ref{E:4}) yields $
\varphi(s)=[z^n]\,\Xi_3(z,s)\sim K(s)\,\frac{4!}{n(n-1)\dots(n-4)}\,
\left(\rho_3(s)^{-1}\right)^n $ with uniform error term and we
accordingly obtain
\begin{equation}\label{E:uu}
f_{X_n}(t)\sim \frac{K(\frac{it}{\sigma_n})}{K(0)}\,
\left[\frac{\rho_3(\frac{it}{\sigma_n})}{\rho_3(0)}\right]^{-n}
e^{-i\frac{\mu_{n}}{\sigma_{n}}t} \ .
\end{equation}
where the error term is uniform for $t$ from any bounded interval.
Taking the logarithm we obtain
\begin{equation}
\ln f_{X_n}(t)\sim\ln\frac{K(\frac{it}{\sigma_n})}{K(0)}- n \,
\ln\frac{\rho_3(\frac{it}{\sigma_n})}{\rho_3(0)} -i \frac{\mu_n}{\sigma_n}t
\ .
\end{equation}
Expanding $g(s)=\ln\frac{\rho_3(s)}{\rho_3(0)}$ in its Taylor series at $s=0$,
(note that $g(0)=0$ holds) yields
\begin{equation}\label{E:kkk}
\ln\frac{\rho_3(\frac{it}{\sigma_{n}})}{\rho_3(0)}=
\frac{\rho_3'(0)}{\rho_3(0)}\frac{it}{\sigma_{n}}-
\left[\frac{\rho_3''(0)}{\rho_3(0)}-
\left(\frac{\rho_3'(0)}{\rho_3(0)}\right)^2\right]
\frac{t^2}{2\sigma^2_{n}}
+{O}(\left(\frac{it}{\sigma_n}\right)^{3})
\end{equation}
and therefore
\begin{equation}\label{E:nun}
\ln f_{n}(t) \sim\ln \frac{K(\frac{it}{\sigma_n})}{K(0)}- n \,
\left\{\frac{\rho_3'(0)}{\rho_3(0)}\frac{it}{\sigma_{n}}-\frac{1}{2}
\left[\frac{\rho_3''(0)}{\rho_3(0)}-
\left(\frac{\rho_3'(0)}{\rho_3(0)}\right)^2\right]
\frac{t^2}{\sigma^2_{n}}
+{O}(\left(\frac{it}{\sigma_n}\right)^{3})\right\}-
\frac{i\mu_{n}t}{\sigma_{n}} \ .
\end{equation}
Claim $2$ implies $\Xi_{3}(z,s)=
{O}\left((\rho_3(s)-z)^4\ln\frac{1}{\rho_3(s)-z}\right)$ is analytic
in $s$ where $s$ is contained in a disc of radius $\epsilon$ around
$0$. Hence $\Xi_{3}(z,s)$ is in particular continuous in $s$ for
$\vert s\vert <\epsilon$ and we can conclude from eq.~(\ref{E:4})
for fixed $t\in ]-\infty, \infty[$
\begin{equation}
\lim_{n\to \infty}\left(\ln K(\frac{it}{\sigma_n})-\ln K(0)\right)=0 \ .
\end{equation}
In view of eq.~(\ref{E:nun}) we introduce
$$
\mu=-\frac{\rho_3'(0)}{\rho_3(0)}, \quad \quad \sigma=\left\{
\left(\frac{\rho_3'(0)}{\rho_3(0)}\right)^2-\frac{\rho_3''(0)}
{\rho_3(0)}\right\}
$$
and eq.~(\ref{E:nun}) becomes
\begin{equation}
\ln f_{X_n}(t)\sim -\frac{\; t^2}{2}
+{O}(\left(\frac{it}{\sigma_n}\right)^{3})
\end{equation}
with uniform error term for $t$ from any bounded interval. This is
equivalent to $\lim_{n\rightarrow\infty}f_{X_n}(t)=
\exp(-\frac{t^2}{2})$ with uniform error term. The
L\'{e}vy-Cram\'{e}r Theorem (Theorem~\ref{T:K}) implies now
eq.~(\ref{E:converge}) and it remains to compute the values for
$\mu$ and $\sigma$ which are given by
\begin{eqnarray}\label{E:mu}
\mu & = & -\frac{\rho_{3}'(0)}{\rho_{3}(0)}=
-\frac{-\frac{3}{2}+\frac{13}{42}
\sqrt{21}}{\frac{5}{2}-\frac{1}{2}\sqrt{21}}=0.39089 \\
\label{E:sigma}
\sigma^2 & = & \mu^2-\frac{\rho_{3}''(0)}{\rho_{3}(0)}=
\mu^2-\frac{1-\frac{94}{441}\sqrt{21}}{\frac{5-\sqrt{21}}{2}}=0.041565
\end{eqnarray}
whence eq.~(\ref{E:concrete0}) and the proof of Theorem~\ref{T:gauss} is
complete.
\end{proof}

\section{The local limit theorem}\label{S:LL}

In this section we complement the central limit theorem presented in the
previous section
$$
\lim_{n\to\infty}\mathbb{P}\left(\frac{X_{n}-\mu_{n}}{\sigma_{n}}<x\right)
=
\frac{1}{\sqrt{2\pi}}\int_{-\infty}^{x}e^{-\frac{t^2}{2}}dt
$$
by considering a ''local'' perspective on the limiting distribution of
$X_{n}$. For the local limit theorem we analyze the difference between
$\mathbb{P}(x \le \frac{X_{n}-\mu_{n}} {\sigma_{n}}<x+1)$ and
$\frac{1}{\sqrt{2\pi}}\int_{x}^{x+1}e^{-\frac{t^2}{2}}dt$ as $n$
tends to infinity.
$X_{n}$ satisfies a local limit theorem on some set $S\subset \mathbb{R}$
if and only if
\begin{equation}\label{E:51}
\lim_{n \to \infty}\sup_{x\in S}\left|\sigma_{n}
\mathbb{P}\left(\frac{X_{n}-\mu_{n}}{\sigma_{n}}=x\right)-\frac{1}{\sqrt{2
\pi}}e^{-\frac{x^2}{2}}\right|=0 \ .
\end{equation}
One key condition formulated in eq.~(\ref{E:wi}) of Theorem~\ref{T:local}
below for proving a local limit theorem is given by
$$
\frac{\varphi_{n}(s)}{\varphi_{n}(0)}\sim \exp(M(s)\beta_{n}+N(s)) \ ,
$$
where $M(s)$ is differentiable and $N(s)$ is continuous in some
$\epsilon$-disc centered at $0$.
In view of eq.~(\ref{E:well}) and eq.~(\ref{E:uu})
in the proof of the central limit theorem, this condition alone implies
the central limit theorem. In other words, the local limit theorem implies
the central limit theorem. We have shown in the introduction that a
central limit theorem does not imply a local limit theorem. Bender observed
in \cite{Bender:73} that the central limit theorem combined with certain
smoothness conditions does imply the local limit theorem. Accordingly, in
order to prove the local limit theorem for $3$-noncrossing RNA structures with
$h$ arcs our strategy will consist in verifying such smoothness conditions
\cite{Hwang:98}.
\begin{theorem}\label{T:local}
Let $\varphi_{n}(s)=\sum_{k}a_{n,k}w^k$ and $w=e^s$. Suppose
\begin{equation}\label{E:wi}
\frac{\varphi_{n}(s)}{\varphi_{n}(0)}\sim \exp(M(s)\beta_{n}+N(s))
\end{equation}
holds uniformly for $\vert s\vert \le \tau$, $s \in \mathbb{C}$ and
$\tau>0$, where the following conditions are satisfied\\
${\sf (i)}$ $M(s)$ is differentiable and $N(s)$ is continuous in
$\vert s\vert<\epsilon$ and furthermore  $M(s)$ and $N(s)$ are
independent of $n$.\\
${\sf (ii)}$ $\beta_{n}$ is independent of $t$, $\beta_{n} \rightarrow \infty$
and $M''(0)>0$;\\
${\sf (iii)}$ there exist constant $\delta$ and $c=c(\delta,r)>0$,
where $0<\delta \le \tau$ such that
\begin{equation}
\left|\frac{\varphi_{n}(r+it)}{\varphi_{n}(r)}\right|=O(\exp(-c\beta_{n}))
\end{equation}
holds uniformly for $-\tau \le r\le \tau$ and $\delta \le \vert t
\vert \le \pi$ as $n$ tends to infinity.\\ Then random variable
$X_{n}$ having distribution $\mathbb{P}(X_{n}=k)=a_{n,k}/a_{n}$
with mean $M'(0)\beta_{n}$ and variance $M''(0)\beta_{n}$
satisfies a local limit theorem on the real set $S=\{x\mid
x=o(\sqrt{\beta_{n}})\}$ i.e.~
\begin{equation}
\lim_{n \to \infty}\sup_{x\in S}\left|\sigma_{n}
\mathbb{P}\left(\frac{X_{n}-\mu_{n}}{\sigma_{n}}=x\right)-\frac{1}{\sqrt{2
\pi}}e^{-\frac{x^2}{2}}\right|=0 \ .
\end{equation}
\end{theorem}

With the help of Theorem~\ref{T:local} we can now prove the local limit
theorem for $3$-noncrossing RNA structures with $h$ arcs.

\begin{theorem}\label{T:local2}
Let ${\sf S}'_3(n,h)$ be the number of $3$-noncrossing RNA structures with
exactly $h$ arcs. Let $X_n$ be the r.v.~having the distribution
\begin{equation}
\forall\; h=0,1,\ldots \lfloor \frac{n}{2}\rfloor,\qquad
\mathbb{P}(X_n=h)=\frac{{\sf S}_3'(n,h)}{{\sf S}_3(n)}
\end{equation}
Then we have for set $S=\{x\mid x=o(\sqrt{n})\}$
\begin{equation}\label{E:52}
\lim_{n \to \infty}\sup_{x\in S}\left|\sqrt{\sigma^2 n}\
\mathbb{P}\left(\frac{X_{n}-n\,\mu}{\sqrt{\sigma^2 n}}
=x\right)-\frac{1}{\sqrt{2
\pi}}e^{-\frac{x^2}{2}}\right|=0  \  ,
\end{equation}
where $\mu=0.39089$ and $\sigma^2  = 0.041565$.
\end{theorem}

\begin{proof}
We will show that $\{\frac{{\sf{S}}'_{3}(n,h)}{{\sf{S}}_{3}(n)}\}$
satisfies the conditions for Theorem \ref{T:local}. For $\vert
s\vert \le \epsilon$, where $\epsilon$ is sufficiently small but
fixed. The crucial equation implying the conditions of
Theorem~\ref{T:local} is eq.~(\ref{E:4}) of the proof of
Theorem~\ref{T:gauss}:
$$
\varphi_{n,3}(s)= K(s)\frac{4!}{n(n-1)\ldots
(n-4)}(\rho_{3}(s)^{-1})^n\left(1-O(\frac{1}{n})\right)\qquad
K(s)\in\mathbb{C} \ ,
$$
holds uniformly for $\vert s\vert<\epsilon$. Therefore we have
\begin{equation}
\frac{\varphi_{n,3}(s)}{\varphi_{n,3}(0)}=
\frac{K(s)}{K(0)}\left(\frac{\rho_{3}(0)}{\rho_{3}(s)}\right)^n
\left(1-O(\frac{1}{n})\right)\sim
\exp\left(n\ln\left(\frac{\rho_{3}(0)}{\rho_{3}(s)}\right)+
\ln\left(\frac{K(s)}{K(0)}\right)\right) \ .
\end{equation}
uniformly for $\vert s\vert< \epsilon$. We set
\begin{equation}
\beta_{n}  =  n,  \quad
M(s)       =  \ln\left(\frac{\rho_{3}(0)}{\rho_{3}(s)}\right)  \
\text{\rm and} \
N(s)       =  \ln\left(\frac{K(s)}{K(0)}\right) \ .
\end{equation}
By construction $t$ is independent of $n$ and clearly $n
\rightarrow \infty$ and $M(s)$ is differentiable and $N(s)$ is
continuous for all $s$ such that $\vert s\vert<\epsilon$.
In addition $M''(0)$ is
analytic for $\vert s\vert<\epsilon$ and we have
$M''(0)=\mu^2-\frac{1-\frac{94}{441}\sqrt{21}}{\frac{5-\sqrt{21}}{2}}=
0.041565>0$. Let $\delta=\epsilon$ and $-\epsilon\le r\le \epsilon$,
we obtain
$$
\frac{\varphi_{n,3}(r+it)}{\varphi_{n,3}(r)}\sim
\exp\left(n\ln\left(\frac{\rho_{3}(r)}
{\rho_{3}(r+it)}\right)+\ln\left(\frac{K(r+it)}{K(r)}\right)\right)
\
$$
uniformly for $-\epsilon \le r\le \epsilon$ and
$\epsilon\le \vert t\vert\le \pi$.
Since $\frac{K(s)}{K(0)}$ yields a constant factor and
$K(s)$ is continuous for $|s|<\epsilon$, it suffices to analyze
$\ln\left(\frac{\rho_{3}(r)}{\rho_{3}(r+it)}\right)$. We observe
$\rho_{3}(s)=\frac{1+4e^{\frac{s}{2}}-\sqrt{12e^s-8e^{\frac{s}{2}}+1}}{2e^s}
\ne 0$ for any complex $s$ where $\vert s\vert< \epsilon$. The singularities
of $\ln(\frac{\rho_{3}(0)}{\rho_{3}(s)})$ correspond to the zeros of
$12e^s-8e^{\frac{s}{2}}+1=(2e^{\frac{s}{2}}+1)
(6e^{\frac{s}{2}}+1)$, that is $e^{\frac{s}{2}}=-\frac{1}{2}$ or
$-\frac{1}{6}$. Observe that for $\vert s \vert < \epsilon$, $\vert
e^{\frac{s}{2}}\vert$ is close to 1. Therefore
$\ln\left(\frac{\rho_{3}(r)}{\rho_{3}(r+it)}\right)$ is analytic for
any $\epsilon\le|t|\le \pi$ and $r \in ]-\epsilon,\epsilon[$ and we
can conclude
$$
\left|\frac{\varphi_{n,3}(r+it)}{\varphi_{n,3}(r)}\right|=O\left|\exp(n\cdot
\ln\left(\frac{\rho_{3}(it)}{\rho_{3}(0)}\right))\right|=
O\left(\exp\left(\mbox{Re}\left(n\cdot
\ln\left(\frac{\rho_{3}(r+it)}{\rho_{3}(r)}\right)\right)\right)\right)
$$
uniformly for $-\epsilon \le r\le \epsilon$ and
$\epsilon\le \vert t\vert\le \pi$.
Taylor expansion of $\ln(\frac{\rho_{3}(r+it)}{\rho_{3}(r)})$ at $0$
shows (see eq.~(\ref{E:kkk})), that the dominant real part of
$\ln(\frac{\rho_{3}(r+it)}{\rho_{3}(r)})$ is given by
$$
\left[\left(\frac{\rho_{3}'(r)}{\rho_{3}(r)}\right)^2-\frac{\rho_{3}''(r)}
{\rho_{3}(r)}\right]\frac{t^2}{2!}<0 \mbox{ for }
r\in]-\epsilon,\epsilon[\ .
$$
Setting $c_{1}=
\left[\frac{\rho_{3}''(r)}{\rho_{3}(r)}-\left(\frac{\rho_{3}'(r)}
{\rho_{3}(r)}\right)^2\right]\frac{\pi^2}{2!}>0$ and
$c_{2}=\left[\frac{\rho_{3}''(r)}{\rho_{3}(r)}-\left(\frac{\rho_{3}'(r)}
{\rho_{3}(r)}\right)^2\right]\frac{\delta^2}{2!}>0$ we can conclude
$$
\left|\frac{\varphi_{n,3}(r+it)}{\varphi_{n,3}(r)}\right|=O(\exp(-c\cdot
n))
$$ for some $0<c_{2}<c<c_{1}$, uniformly for $-\epsilon \le r\le \epsilon$ and
$\epsilon\le \vert t\vert\le \pi$ and Theorem~\ref{T:local} applies, whence
Theorem~\ref{T:local2}.
\end{proof}

\section{Appendix}


{\bf Proof of Lemma~\ref{L:func}.} First we observe that for
$x,w\in [-1,1]$ the term $w^2x^2-x+1$ is strictly positive. We set
\begin{equation}
F_k(x,w)=\sum_{n\ge 0} \sum_{h\le n/2}{\sf S}_k'(n,h)w^{2h}x^n
\end{equation}
and compute
\begin{eqnarray*}
F_k(x,w) & = & \sum_{n\ge 0}\sum_{h\le
n/2}\sum_{j=0}^h(-1)^j\binom{n-j}{j}
\binom{n-2j}{2(h-j)} f_k(2(h-j),0)w^{2h}x^n \\
& = & \sum_{n\ge 0}\sum_{j\le n/2}\sum_{h=j}^{n/2}(-1)^j\binom{n-j}{j}
\binom{n-2j}{2(h-j)}f_k(2(h-j),0)w^{2h}x^n \\
& = & \sum_{j\ge 0}\sum_{n\ge
2j}\sum_{h=j}^{n/2}(-1)^j\binom{n-j}{j}
\binom{n-2j}{2(h-j)}f_k(2(h-j),0)w^{2h}x^n \\
& = & \sum_{j\ge 0}(-1)^j\frac{(wx)^{2j}}{j!}\sum_{n\ge 2j}(n-j)!
\sum_{h=j}^{n/2}\binom{n-2j}{2(h-j)}f_k(2(h-j),0)\frac{w^{2(h-j)}}{(n-2j)!}
x^{n-2j} \ .
\end{eqnarray*}
We shift summation indices $n'=n-2j$ and $h'=h-j$ and derive for the rhs
the following expression
\begin{eqnarray*}
& = & \sum_{j\ge 0}(-1)^j\frac{(wx)^{2j}}{j!}\sum_{n'\ge 0}(n'+j)!
\sum_{h=j}^{n/2}\binom{n'}{2(h-j)}f_k(2(h-j),0)\frac{w^{2(h-j)}}{n'!}
x^{n-2j} \\
& = & \sum_{j\ge 0}(-1)^j\frac{(wx)^{2j}}{j!}\sum_{n'\ge 0}(n'+j)!
\left\{
\sum_{h'=0}^{n/2-j=n'/2}\binom{n'}{2h'}f_k(2h',0)w^{2h'}\right\}
\frac{x^{n'}}{n'!}
\end{eqnarray*}
The idea is now to interpret the term
$\sum_{h'=0}^{n'/2}\binom{n'}
{2h'}f_k(2h',0)w^{2h'}\frac{x^n}{n!}$ as a product of the two
power series $e^x$ and $\sum_{n\ge
0}f_k(2n,0)\frac{(wx)^{2n}}{(2n)!}$:
\begin{eqnarray*}
\sum_{\ell\ge 0}\frac{x^\ell}{\ell!}\sum_{n\ge
0}f_k(2n,0)\frac{(wx)^{2n}} {(2n)!} & = & \sum_{n'\ge 0}
\sum_{2n+\ell=n'} \left\{\frac{1}{\ell!}\frac{1}{(2n)!}
f_k(2n,0)w^{2n}\right\}x^{n'}\\
& = & \sum_{n'\ge 0}\left\{\sum_{n=0}^{n'/2}
\binom{n'}{2n}f_k(2n,0)w^{2n}\right\}\frac{x^{n'}}{n'!} \ .
\end{eqnarray*}
We set
$\eta_{n'}=\left\{\sum_{h'=0}^{n'/2}\binom{n'}{2h'}f_k(2h',0)
w^{2h'}\right\} $. By assumption we have $\vert x\vert<\rho_k(w)$
and we next derive, using the Laplace transformation and
interchanging integration and summation
\begin{equation}\label{E:on1}
\sum_{n'\ge 0}(n'+j)!\eta_n \frac{x^{n'}}{n'!} = \int_{0}^{\infty}
\sum_{n'\ge 0}\eta_{n'} \frac{(xt)^{n'}}{n'!} t^je^{-t} dt \ .
\end{equation}
Since $\vert x\vert<\rho_k(w)$ the above transformation is valid
and using
\begin{equation}
 \sum_{n'\ge 0}\left\{\sum_{n=0}^{n'/2}
\binom{n'}{2n}f_k(2n,0)w^{2n}\right\}\frac{x^{n'}}{n'!}=
\sum_{\ell\ge 0}\frac{x^\ell}{\ell!}\sum_{n\ge
0}f_k(2n,0)\frac{(wx)^{2n}} {(2n)!}
\end{equation}
we accordingly obtain
\begin{eqnarray}
\sum_{n'\ge 0}\eta_{n'} \frac{(xt)^{n'}}{n'!} t^je^{-t} dt & = &
 \int_{0}^{\infty}e^{tx}\sum_{n\ge 0}f_k(2n,0)\frac{(wxt)^{2n}}{(2n)!}
t^je^{-t}dt \ .
\end{eqnarray}
The next step is to substitute the term $\sum_{n'\ge
0}(n'+j)!\eta_n \frac{x^{n'}}{n'!}$ in eq.~(\ref{E:on1}), whence
consequently
\begin{eqnarray*}
F_k(x,w) & = &  \sum_{j\ge 0}(-1)^j\frac{(wx)^{2j}}{j!}
 \int_{0}^{\infty}e^{tx}\sum_{n\ge 0}f_k(2n,0)\frac{(wxt)^{2n}}
{(2n)!}t^je^{-t}dt\\
& = &  \int_{0}^{\infty}\sum_{j\ge 0}(-1)^j\frac{(wx)^{2j}}{j!}
e^{tz}\sum_{n\ge 0}f_k(2n,0)\frac{(wxt)^{2n}} {(2n)!}t^je^{-t}dt \
.
\end{eqnarray*}
The summation over the index $j$ is just an exponential function and we
derive
\begin{eqnarray*}
& = &\int_{0}^{\infty} e^{-(w^2x^2-x+1)t} \sum_{n\ge
0}f_k(2n,0)\frac{(wxt)^{2n}}
{(2n)!}dt \\
& = &\int_{0}^{\infty} e^{-(w^2x^2-x+1)t} \sum_{n\ge
0}f_k(2n,0)\frac{1}{(2n)!} \left(\frac{wx}{w^2x^2-x+1}\right)^{2n}
((w^2x^2-x+1)t)^{2n}dt
\end{eqnarray*}
We proceed by transforming the integral introducing
$u=(w^2x^2-x+1)t$, i.e.~ $dt=(w^2x^2-x+1)^{-1}du$ and accordingly
arrive at
\begin{eqnarray*}
F_k(x,w) & = &  \sum_{n\ge 0}f_k(2n,0)\frac{1}{(2n)!}
\left(\frac{wx}{w^2x^2-x+1}\right)^{2n} \int_{0}^{\infty}
e^{-(w^2x^2-x+1)t}  ((w^2x^2-x+1)t)^{2n}
dt \\
& = &  \sum_{n\ge 0}f_k(2n,0)\frac{1}{(2n)!}
\left(\frac{wx}{w^2x^2-x+1}\right)^{2n}
\frac{1}{w^2x^2-x+1} (2n)! \\
& = & \frac{1}{w^2x^2-x+1}  \sum_{n\ge 0}f_k(2n,0)
\left(\frac{wx}{w^2x^2-x+1}\right)^{2n} \ ,
\end{eqnarray*}
In particular for $w=1$
\begin{eqnarray}\label{E:L1}
\sum_{n \ge 0}{\sf S}_{k}(n)x^n=\frac{1}{x^2-x+1}\sum_{n \ge
0}f_{k}(2n,0)\left(\frac{x}{x^2-x+1}\right)^{2n}
\end{eqnarray}
holds for any $x\in\mathbb{R}$, satisfying $\vert
x\vert<\rho_{k}(1)$, and where $\rho_{k}(1)$ is the radius of
convergence of the power series $\sum_{n \ge 0}{\sf S}_{k}(n)z^n$
over $\mathbb{C}$, that is eq.~(\ref{E:L1}) holds for $x\in
]-\rho_{k}(1),\rho_{k}(1)[$ . From complex analysis we know that
any two functions that are analytic at $0$ and coincide on an open
interval which includes $0$ are identical. Therefore
eq.~(\ref{E:L1}) holds for $z\in\mathbb{C}$, $\vert
z\vert<\rho_k(1)$, and the proof of the lemma is complete.
$\square$

{\bf Acknowledgments.}
We are grateful to Prof.~Jason Gao for helpful discussions.
This work was supported by the 973 Project, the PCSIRT Project of the
Ministry of Education, the Ministry of Science and Technology, and
the National Science Foundation of China.

\bibliography{central}
\bibliographystyle{plain}


\end{document}